\documentclass[11pt,a4paper]{amsart}
\usepackage[a4paper,margin=26mm]{geometry}
\usepackage{etoolbox}
\usepackage{amssymb, color, graphicx, comment}
\usepackage[hidelinks]{hyperref}
\numberwithin{equation}{section}
\newtheorem{theorem}{Theorem}[section]
\newtheorem{lemma}[theorem]{Lemma}

\newtheorem{corollary}[theorem]{Corollary}
\theoremstyle{definition}
\newtheorem{remark}[theorem]{Remark}
\newtheorem{definition}{Definition}[section]
\theoremstyle{remark}
\newcommand{\Tr}{\operatorname{Tr}}
\newcommand{\E}{\mathcal E}
\newcommand{\ind}{\mathbf 1}
\newcommand{\BoxC}[2]{[#1,#2]_{\mathrm{Box}}}
\newcommand{\CE}[2]{[#1,#2]_{\mathrm{CE}}}
\newcommand{\HC}[2]{[#1,#2]_{\mathrm{HC}}}

\newcommand{\Hardy}{\mathbf I}
\newcommand{\Cartprod}{%
  \mathop{\vcenter{\hbox{\LARGE$\times$}}}\displaylimits
}
\title{Weighted Carleson embedding on product of trees in general dimension}
\author{}
\makeatletter
\renewenvironment{proof}[1][\proofname]{%
  \par\pushQED{\qed}%
  \normalfont
  \topsep6\p@\@plus6\p@\relax
  \trivlist
  \item[\hskip\labelsep\bfseries #1\@addpunct{.}]\ignorespaces
}{%
  \popQED\endtrivlist\@endpefalse
}
\makeatother

\author[M. Paparizos]{Michail Paparizos}
\address{Department of Mathematics, Michigan State University, East Lansing, Michigan, USA}
\email{paparizo@msu.edu}

\author[G. Psaromiligkos]{Georgios Psaromiligkos}
\address{Department of Mathematics, University of Thessaly, Lamia, Greece}
\email{gpsaromiligkos@uth.gr}
\begin{document}

\begin{abstract}
We prove the equivalence of the Carleson embedding on $n$-trees (arbitrary Cartesian product of trees) and the box condition, for arbitrary nonnegative measures and generalized product weights. For product weights, we obtain the sharp coefficient $4^n$. This extends the bi-tree and tri-tree embedding theorems of \cite{AMPVZ} and \cite{MPVZ} to every dimension. We also give two applications of the Carleson embedding: we characterize Carleson measures for the classical Dirichlet space on the polydisc $\mathbb D^n$ and deduce criteria for bounded composition operators from this space into weighted Bergman spaces, generalizing the results from \cite{MPVZ} and \cite{Beslikas2025}. Finally, we define a generalized logarithmic Dirichlet space and characterize its Carleson measures.
\end{abstract}

\maketitle

\section{Introduction and main results}\label{sec:main}

Is the box condition on multi-trees (Cartesian products of $n$ trees) sufficient for the Carleson embedding? This is a question that was left unanswered in \cite{AMPVZ} and \cite{MPVZ}. One should note that this embedding is different from that of Sun-Yung ~A.~Chang~\cite{Chang1979}, since here the measure on the multi-tree is arbitrary, with the weight having the product form.

More precisely, let $T^n=T_1\times\cdots\times T_n$ be a product of finite rooted trees, ordered coordinatewise, so that $\beta\le\alpha$ means that $\beta$ is a descendant of $\alpha$. For a non-negative vertex measure $\mu$ and a non-negative weight $w$, we consider the embedding
\begin{gather*}
\sum_{\alpha\in T^n}w(\alpha)
\bigg|\sum_{\beta\le\alpha}f(\beta)\mu(\beta)\bigg|^2
\le C\sum_{\beta\in T^n}|f(\beta)|^2\mu(\beta).
\end{gather*}
Testing with the characteristic function of a descendant set, we get the necessary box condition:
\begin{gather*}
\sum_{\alpha\le\gamma}w(\alpha)
\bigg(\sum_{\beta\le\alpha}\mu(\beta)\bigg)^2
\le B\sum_{\beta\le\gamma}\mu(\beta)
\qquad(\gamma\in T^n).
\end{gather*}

In Theorem~\ref{thm:product} we prove that, whenever
$w(\alpha)=\prod_{j=1}^n w_j(\alpha_j)$ and the measure $\mu$ is arbitrary, this testing condition implies the embedding with $C=4^nB$. We also extend this result, in Theorem \ref{thm:corner}, to a more general class of weights, called corner weights.

The one-parameter and multiparameter settings behave quite differently. In one parameter, testing on descendant sets is a sufficient condition for the classical Carleson embedding theorem and its weighted forms, as it can be seen for example in 
\cite{ARS2002},\cite{Carleson1962} and \cite{PsaromiligkosThesis}. On a single tree, descendant sets are either disjoint or nested, and this makes the standard maximal-function argument work. For products of trees, this is no longer the case, since two descendant rectangles may intersect without one containing the other. In such cases, the corresponding maximal operator may not be bounded with respect to an arbitrary measure, and the standard one-parameter argument does not work.

In the dyadic setting, take the underlying measure to be product Lebesgue measure $m_n$ on $[0,1]^n$ and allow arbitrary vertex coefficients $w(R)$. Writing $c_R=w(R)m_n(R)^2$, the box condition becomes $\sum\limits_{R\subseteq R_0}c_R\lesssim m_n(R_0)$ for every dyadic rectangle $R_0$. Carleson's counterexample \cite{Carleson1974} shows that this condition is insufficient even in two parameters. Sun-Yung A.~Chang \cite{Chang1979} obtained the appropriate stronger condition by testing unions of rectangles:
$\sum\limits_{R\subseteq\Omega}c_R\lesssim m_n(\Omega),
$ for finite unions $\Omega$ of dyadic rectangles. The associated embedding belongs to the theory of separately harmonic Poisson extensions. Chang and Fefferman \cite{ChangFefferman1980} developed the corresponding product $H^1$--$\mathrm{BMO}$ duality, in which testing over such unions is essential.

In this setting, the measure is a product measure and the coefficient weight is unrestricted. The bi-tree and tri-tree theorems address the complementary problem. that is, the measure $\mu$ is arbitrary, while $w$ has product structure. Carleson's example uses a nonproduct coefficient weight and therefore lies outside that class. At the same time, counterexamples for general pairs $(w,\mu)$ show that some structure is needed: even testing over unions of rectangles need not characterize the embedding when both inputs are unrestricted \cite{MPVcounterexamples}. This motivates our search for conditions on the weight under which the box condition remains sufficient for arbitrary measures \(\mu\).

This choice of assumptions is dictated by applications in several complex variables. A nonnegative measure $\nu$ is a Carleson measure for a Hilbert space $\mathcal H$ of holomorphic functions if
$$
\int |f|^2\,d\nu\le C_\nu\|f\|_{\mathcal H}^2
\qquad(f\in\mathcal H).
$$

Such measures connect embedding theorems with interpolation, multipliers, and capacity. In one complex variable, the results of Carleson \cite{Carleson1962} and Stegenga \cite{Stegenga1980} provide the basic examples, while tree models give a discrete approach to these questions \cite{ARS2002,ARSW2011}. For the Dirichlet space on the bidisc, Arcozzi, Mozolyako, Perfekt, and Sarfatti \cite{AMPS2023} obtained a characterization through bi-parameter Bessel capacity.

For the Dirichlet spaces treated in \cite{AMPVZ,MPVZ}, the product structure of the reproducing kernel leads to a product weight in the dyadic model. The analytic measure itself is arbitrary, and its masses on product Whitney cells produce an arbitrary vertex measure. In particular, the classical Dirichlet space corresponds to $w\equiv1$. Thus, the natural discrete problem arising from complex analysis has a product coefficient weight and an unrestricted measure. The analytic-to-dyadic reduction of \cite[Theorem~2.2]{MPVZ} is already available in every dimension for the classical Dirichlet space. The remaining higher-dimensional difficulty was the discrete embedding theorem, which we obtained in this paper. 

The earlier proofs on the bi-tree and tri-tree used potential theory and combinatorial energy estimates. One obstacle is that potentials on product trees do not satisfy the classical maximum principle. The arguments in \cite{AMPVZ} and \cite{MPVZ} overcome this difficulty using small-energy majorization and related surrogate maximum principles. These estimates suffice  to connect box testing, Carleson conditions, restricted energy conditions, and the full embedding in dimensions two and three, but the additional mixed interactions that appear with four or more parameters prevent the same argument from going through. In fact, the counterexamples in \cite{MPVfailure} show that certain extensions of the auxiliary majorization estimates cannot hold. See also \cite[Section~8]{MPVZ}.
However, the embedding problem in higher dimensions remained open.

The proofs of Theorems~\ref{thm:product} and~\ref{thm:corner} follow the operator formulation used in used in Psaromili-gkos' thesis
\cite[Section~3.1]{PsaromiligkosThesis}, where a suggestion of
M.~Christ led to a one-dimensional proof based on Schur's test.
Composing the embedding operator with its adjoint produces a
positive operator, whose norm is controlled by constructing
an auxiliary function. In higher dimensions, however, this construction does not extend
\cite[Section~7.1]{PsaromiligkosThesis}. Our approach is based on estimates for traces of powers of certain matrices, using the product structure
of the weight. This gives the sharp coefficient $4^n$ for product
weights in every dimension. The same argument extends to nonproduct weights satisfying a \textit{corner}
condition. Under the unit corner condition, the coefficient $4^n$
remains optimal among nonproduct weights when $n\ge2$.

In Theorem~\ref{thm:capacity-characterization}, we use these embedding estimates to characterize a discrete Hardy capacity for corner weights in every dimension. Earlier work covers the unweighted tri-tree case (e.g.
\cite[Chapter~4]{PsaromiligkosThesis}) and related capacity estimates
(e.g \cite[Section~4.3.2]{MPVZ}). 

The analytic-to-dyadic reduction provides a  characterization of Carleson measures for $\mathcal D(\mathbb D^n)$ in terms of uniform box tests over rotated coordinate lattices, for arbitrary finite nonnegative Borel measures. Applying this characterization to the push-forward measure of weighted Bergman volume, these tests characterize bounded composition operators from $\mathcal D(\mathbb D^n)$ into weighted Bergman spaces, extending Beslikas' bidisc result \cite{Beslikas2025} to every dimension. 

Finally, we define a logarithmic Dirichlet space, where we apply the box condition to characterize the Carleson measures. 

\subsection*{Trees, measures, and the box condition}\label{subsec:model}

Let $T$ be a nonempty finite rooted tree and let $\alpha\in T$ be a vertex. We order $T$ by
\[
\alpha\le\beta
\quad\Longleftrightarrow\quad
\beta\text{ lies on the path from }\alpha\text{ to the root}.
\]

The root is the maximal vertex. A \emph{leaf} is a vertex with no children, equivalently a minimal vertex. The boundary $\partial T$ is the set of all leaves.

For \(n\ge1\), an \emph{$n$-tree} is the Cartesian product
\[
T^n=T_1\times\cdots\times T_n,
\]
with the coordinatewise order
\[
\alpha\le\beta
\quad\Longleftrightarrow\quad
\alpha_j\le\beta_j\text{ in }T_j\text{ for every }j\in\{1,\ldots,n\}.
\]

The boundary of the $n$-tree is
\[
\partial T^n=\partial T_1\times\cdots\times\partial T_n.
\]

For $\alpha\in T^n$ define
\begin{gather*}
\partial\alpha=\{\omega\in\partial T^n:\omega\le\alpha\}.
\end{gather*}

For a measure $\mu$ on the $n$-tree, we write $\mu(\beta)$ for the mass at the vertex $\beta$. The Hardy operator and its adjoint are given by
\begin{gather*}
(\Hardy h)(\alpha)=\sum_{\beta\ge\alpha}h(\beta)\ \ \text{and}\ \ 
(\Hardy^*h)(\alpha)=\sum_{\beta\le\alpha}h(\beta).
\end{gather*}

In particular, for each vertex $\alpha$, the total mass of $\alpha$ and its descendants, that is, the cumulative mass of the measure $\mu$ is
\begin{gather*}
(\Hardy^*\mu)(\alpha)
=\sum_{\beta\le\alpha}\mu(\beta)
\end{gather*}

If $\mu$ is supported on the boundary, this equals $\mu(\partial\alpha)$. For a complex measure $\sigma$ on $T^n$, its weighted energy is
\begin{gather*}
\E_w[\sigma]=\sum_{\alpha\in T^n}w(\alpha)
  |(\Hardy^*\sigma)(\alpha)|^2.
\end{gather*}

Thus for a finite nonnegative measure $\mu$ on $T^n$, and $w:T^n\to[0,\infty)$ a weight: 
\begin{gather*}
\E_w[f\mu]
=\sum_{\alpha\in T^n}w(\alpha)
\left|\sum_{\beta\le\alpha}f(\beta)\mu(\beta)\right|^2.
\end{gather*}

Let $\BoxC{w}{\mu}$ be the least nonnegative constant of the \emph{box condition}:
\begin{equation}\label{eq:boxconstant}
\sum_{\alpha\le\gamma}w(\alpha)(\Hardy^*\mu)(\alpha)^2
\le\BoxC{w}{\mu}\,(\Hardy^*\mu)(\gamma)
\qquad(\gamma\in T^n)
\end{equation}
and $\CE{w}{\mu}$ be the least nonnegative constant for the \emph{Carleson embedding inequality}:
\begin{equation}\label{eq:CEconstant}
\E_w[f\mu]\le\CE{w}{\mu}\,\|f\|_{L^2(\mu)}^2
\qquad(f\in L^2(\mu)).
\end{equation}

As mentioned before, the embedding \eqref{eq:CEconstant} implies the box condition \eqref{eq:boxconstant}. In fact, fix $\gamma$ and set $f(\beta)=\ind_{\{\beta\le\gamma\}}$. Then, $\|f\|_2^2=(\Hardy^*\mu)(\gamma)$. Hence $\BoxC{w}{\mu}\le\CE{w}{\mu}$ for every weight. 

In this paper, we answer the reverse estimate, uniformly in the depth and size of the $n$-trees. In particular, we prove 
\begin{theorem}\label{thm:product}
Let $n\ge1$, let $\mu$ be a finite nonnegative measure on $T^n$, and suppose
\[
w(\alpha)=\prod_{j=1}^n w_j(\alpha_j),
\]
where \(w_j:T_j\to[0,\infty)\) are weights. Then
\begin{gather*}
\CE{w}{\mu}\le4^n\BoxC{w}{\mu}.
\end{gather*}
The coefficient $4^n$ is sharp over these weights.
\end{theorem}

We will deduce it from the more general Theorem~\ref{thm:corner}, which will be stated once we introduce the necessary terminology.

For $\alpha,\beta\in T^n$, write $\alpha\sim\beta$ if for every $j\in\{1,\ldots,n\}$,
\[
\alpha_j\le\beta_j\ \text{or}\ \beta_j\le\alpha_j.
\]
Equivalently, $\partial\alpha\cap\partial\beta\ne\varnothing$.
For such pairs define
\begin{gather*}
(\alpha\wedge\beta)_j=
\begin{cases}\alpha_j,&\alpha_j\le\beta_j,\\\beta_j,&\beta_j\le\alpha_j,\end{cases}
\qquad
(\alpha\vee\beta)_j=
\begin{cases}\beta_j,&\alpha_j\le\beta_j,\\\alpha_j,&\beta_j\le\alpha_j.\end{cases}
\end{gather*}
The definitions agree when the coordinate vertices are equal. Note that
\[
\{\gamma:\gamma\le\alpha,\ \gamma\le\beta\}
=\{\gamma:\gamma\le\alpha\wedge\beta\}.
\]

A \emph{coordinate split} is a subset $S\subseteq\{1,\ldots,n\}$, with complement $S^c=\{1,\ldots,n\}\setminus S$. Set
{
\[
\mathcal T_S=\Cartprod_{j\in S}T_j,\qquad\mathcal T_{S^c}=\Cartprod_{j\in S^c}T_j. 
\]

We index each block by its original tree coordinates. Thus, if $a=(a_j)_{j\in S}$ and $v=(v_j)_{j\in S^c}$, define $(a,v)\in T^n$ by
\begin{gather}
(a,v)_j=\begin{cases}a_j,&j\in S,\\v_j,&j\in S^c.\end{cases}
\label{eq:block-assembly}
\end{gather}

In case $S=\varnothing$ or $S=\{1,\ldots,n\}$ then $(a,v)=v$ and $(a,v)=a$ respectively.

\begin{definition}[Corner condition]\label{def:corner}
A weight $w\ge0$ satisfies the corner condition for the split $S$ with constant $K_S\ge1$ if, for all $a,u\in\mathcal T_S$ and $b,v\in\mathcal T_{S^c}$,
\begin{gather}
w(a,v)w(u,b)\le K_S\ w(a,b)w(u,v)
\quad(a\le u,\ b\le v).
\label{eq:corner}
\end{gather}
or, equivalently,
\[w(\alpha)w(\beta)
\le K_S\ w(\alpha\wedge\beta)w(\alpha\vee\beta)
\qquad(\alpha\sim\beta).
\] 
\end{definition}

For $S=\varnothing$ and $S=\{1,\ldots,n\}$ the inequality is equality, so we take $K_\varnothing=K_{\{1,\ldots,n\}}=1$.
Complementary splits have the same optimal constant.

We say that $w$ satisfies the \emph{unit corner condition} if
$K_S\equiv1$, that is,
\begin{gather}
w(\alpha)w(\beta)
\le w(\alpha\wedge\beta)w(\alpha\vee\beta)
\qquad(\alpha\sim\beta).
\label{eq:unitcorner}
\end{gather}

This condition is more general than the product weight assumption
in the earlier bi-tree and tri-tree theorems of \cite{AMPVZ} and \cite{MPVZ},
since product weights satisfy it with equality. Our main result is the following:

\begin{theorem}\label{thm:corner}
Let $\mu$ be a finite nonnegative measure on $T^n$, $n\geq 1$. If $w$ satisfies \eqref{eq:corner} for every coordinate split with finite constants $K_S$, then
\begin{gather}
\CE{w}{\mu}\le\kappa(w)\BoxC{w}{\mu},
\label{eq:directcorner}
\end{gather}
where \(\displaystyle \kappa(w)=\left(\sum_{S\subseteq\{1,\ldots,n\}}\sqrt{K_S}\right)^2. \)
If the split constants satisfy $K_S\le K$ for every nonempty proper subset $S\subset\{1,\ldots,n\}$, with a common bound $K\ge1$, then
\begin{gather}
\CE{w}{\mu}
\le\bigl(2+(2^n-2)\sqrt K\bigr)^2\BoxC{w}{\mu}
\le4^nK\BoxC{w}{\mu}.
\label{eq:globalK}
\end{gather}
In particular, the unit corner condition \eqref{eq:unitcorner}
gives the sharp coefficient $4^n$.
\end{theorem}

Before we turn to the proofs of these two theorems, let us give some applications to the Dirichlet and weighted Bergman spaces on the polydisc.
First, Theorem \ref{thm:polydisc-dirichlet} characterizes Carleson measures for the classical Dirichlet space in every dimension, extending the already known results for two and three dimensions. 

Write $\mathbb D=\{z\in\mathbb C:|z|<1\}$ and
$\mathbb T=\mathbb R/\mathbb Z$, with the identification
$t\mapsto e^{2\pi it}$. For a positive integer $n$, define
$\mathcal D(\mathbb D^n)$ to be the space of holomorphic functions
$f(z)=\sum\limits_{k\in\mathbb N_0^n}a_kz^k$ such that
\begin{gather}
\|f\|_{\mathcal D(\mathbb D^n)}^2
=\sum_{k\in\mathbb N_0^n}|a_k|^2
  \prod_{j=1}^n(k_j+1)<\infty.
\label{eq:polydisc-dirichlet-norm}
\end{gather}
For a finite nonnegative Borel measure $\nu$ on the open polydisc,
let $[\nu]_{\mathrm{CM}}\in[0,\infty]$ denote its Carleson measure norm, that is, the least constant in
\begin{gather*}
\int_{\mathbb D^n}|f|^2\,d\nu
\le [\nu]_{\mathrm{CM}}\|f\|_{\mathcal D(\mathbb D^n)}^2
\qquad(f\in\mathcal D(\mathbb D^n)).
\end{gather*}
The measure is a Carleson measure for $\mathcal D(\mathbb D^n)$
when this constant is finite.

For $\theta\in\mathbb T$, let $\mathcal D^\theta$ consist of the
half-open arcs
\begin{gather*}
I_{m,k}^\theta
=\{\theta+t\!\!\!\! \pmod 1:k2^{-m}\le t<(k+1)2^{-m}\},
\qquad m\ge0,\quad 0\le k<2^m.
\end{gather*}
The root is the full circle, the arc length is normalized so that
$|\mathbb T|=1$ and the descendants are dyadic subarcs. For $I\in\mathcal D^\theta$, let
\begin{equation}\label{Stheta}
S_\theta(I)=\{re^{2\pi it}:t\in I,\ 1-|I|\le r<1\}.
\end{equation}
For $\omega=(\theta_1,\ldots,\theta_n)\in\mathbb T^n$, set
$\mathcal R^\omega=\mathcal D^{\theta_1}\times\cdots\times\mathcal D^{\theta_n}$ and define
\[
S_\omega(\alpha)=\Cartprod_{j=1}^n S_{\theta_j}(I_j)     
\qquad\bigl(\alpha=I_1\times\cdots\times I_n\in\mathcal R^\omega\bigr).
\]

Let $[\nu]_{\mathrm{Box}}$ be the least nonnegative constant for which
\begin{gather}
\sum_{\substack{\beta\in\mathcal R^\omega\\\beta\le\alpha}}
\nu(S_\omega(\beta))^2\le[\nu]_{\mathrm{Box}}\nu(S_\omega(\alpha)),
\label{eq:polydisc-uniform-box-test}
\end{gather}
for \(\omega\in\mathbb T^n,\text{ and } \alpha\in\mathcal R^\omega\). If there is not such a finite constant, we set $[\nu]_{\mathrm{Box}}=+\infty$. In this setting \(\beta\le \alpha\) means coordinate dyadic inclusion.

\begin{theorem}\label{thm:polydisc-dirichlet}
Let $n\ge1$, and let $\nu$ be a finite nonnegative Borel measure
on $\mathbb D^n$. Then $\nu$ is a Carleson measure for
$\mathcal D(\mathbb D^n)$ if and only if
$[\nu]_{\mathrm{Box}}<\infty$. More precisely,
\begin{gather}
\frac{1}{b_n}[\nu]_{\mathrm{Box}}
\le [\nu]_{\mathrm{CM}}
\le a_n4^n[\nu]_{\mathrm{Box}},
\label{eq:polydisc-main-bounds}
\end{gather}
where $a_n,b_n>0$ depend only on $n$.
\end{theorem}

Our second application gives a Carleson box criterion
for the composition problem involving weighted Bergman spaces in every dimension, extending the bidisc result of Beslikas in \cite{Beslikas2025}. 

More precisely, let $\beta>-1$ and set
\begin{gather*}
dV_\beta(z)=(\beta+1)^n\prod_{j=1}^n(1-|z_j|^2)^\beta\,dV(z),
\end{gather*}
where $dV$ is normalized volume measure on $\mathbb D^n$.
Let the space $A^2_\beta(\mathbb D^n)$ consist of the holomorphic
functions in $L^2(V_\beta)$. Note that for $\beta=0$, we have the usual Bergman
space. For a holomorphic function $\Phi:\mathbb D^n\to\mathbb D^n$,
we set $\nu_\Phi=V_\beta\circ\Phi^{-1}=\Phi_*V_\beta$, that is,
$\nu_\Phi(E)=V_\beta(\Phi^{-1}(E))$ for Borel sets $E$. Beslikas obtained the bidisc criterion by applying the capacity characterization in \cite{AMPS2023} to the same measure
$V_\beta\circ\Phi^{-1}$. 

The next corollary allows each coordinate function of \(\Phi\) to depend on all variables.

\begin{corollary}\label{cor:polydisc-composition}
The composition operator
$C_\Phi f=f\circ\Phi$ is bounded from
$\mathcal D(\mathbb D^n)$ to $A^2_\beta(\mathbb D^n)$ if and only if
$[\nu_\Phi]_{\mathrm{Box}}<\infty$. Its squared norm satisfies
\begin{gather*}
\frac{1}{b_n}[\nu_\Phi]_{\mathrm{Box}}
\le\|C_\Phi\|^2
\le a_n4^n[\nu_\Phi]_{\mathrm{Box}}.
\end{gather*}
\end{corollary}

The proofs of our main result, Theorem \ref{thm:corner}, and its corollary, Theorem \ref{thm:product}, are given in Section \ref{sec:proofs}. Section \ref{sec:capacity} is devoted to the proof of Theorem \ref{thm:polydisc-dirichlet} and Corollary \ref{cor:polydisc-composition}. In Section \ref{sec:capacity}, we also develop the capacitary relations with the embedding inequalities.

\section{Proof of the Main Result}\label{sec:proofs}

We work over complex Hilbert spaces, with inner products linear in the first argument. Define $\mathcal A:L^2(\mu)\to\ell^2(T^n)$ by
\begin{gather*}
(\mathcal Af)(\alpha)
=\sqrt{w(\alpha)}(\Hardy^*(f\mu))(\alpha)
=\sqrt{w(\alpha)}\sum_{\beta\le\alpha}f(\beta)\mu(\beta).
\end{gather*}
Recalling that $\langle\mathcal Af,h\rangle=\langle f,\mathcal A^{*}h\rangle$ and having changed the order of the summation in $\langle\mathcal Af,h\rangle$, we get the formula for the adjoint operator $\mathcal  A^*$:
\begin{gather*}
(\mathcal A^*h)(\beta)=\sum_{\alpha\ge\beta}\sqrt{w(\alpha)}h(\alpha).
\end{gather*}
Let $|T^n|=\#\{\text{total vertices of }T^n\}$. The Hermitian matrix $G=\mathcal{A}\mathcal{A}^{*}$ is a $|T^n|\times|T^n|$ matrix with nonnegative entries 
\begin{gather*}
G(\alpha,\beta)=\sqrt{w(\alpha)w(\beta)}
\sum_{\substack{\gamma\le\alpha\\\gamma\le\beta}}\mu(\gamma)
=\begin{cases}
\sqrt{w(\alpha)w(\beta)}\,(\Hardy^*\mu)(\alpha\wedge\beta),&\alpha\sim\beta,\\
0,&\alpha\not\sim\beta.
\end{cases}
\end{gather*}
Note also that $G$ is positive semidefinite because
$\langle Gh,h\rangle=\|\mathcal A^*h\|_{L^2(\mu)}^2\ge0$. Consequently,
\begin{gather}
\|G\|=\|\mathcal A\|^2=\CE{w}{\mu}.
\label{eq:norms}
\end{gather}

Now, we define the matrix $F$ with entries 
$$
F(\alpha,\beta)
=G(\alpha,\beta)\mathbf 1_{\{\beta\le\alpha\}}
=\sqrt{w(\alpha)w(\beta)}\,(\Hardy^*\mu)(\beta)\mathbf 1_{\{\beta\le\alpha\}}.
$$ 

\begin{lemma}\label{lem:energy}
For every nonnegative weight $w$,
\begin{gather*}
\|F\|^2\le\BoxC{w}{\mu}\,\|G\|.
\end{gather*}
\end{lemma}

\begin{proof}
Matrix multiplication gives
\begin{gather*}
(FF^*)(\alpha,\beta)=\sqrt{w(\alpha)w(\beta)}
\sum_{\substack{\gamma\le\alpha\\\gamma\le\beta}}w(\gamma)(\Hardy^*\mu)(\gamma)^2.
\end{gather*}
If $\alpha\sim\beta$, their common descendants are exactly the descendants of $\alpha\wedge\beta$. By inequality \eqref{eq:boxconstant},
\begin{equation}\label{FFstar}
(FF^*)(\alpha,\beta)
\le\BoxC{w}{\mu}\sqrt{w(\alpha)w(\beta)}\,(\Hardy^*\mu)(\alpha\wedge\beta)
=\BoxC{w}{\mu}\,G(\alpha,\beta).
\end{equation}
If $\alpha\not\sim\beta$, at least one coordinate pair has no common descendant, so $FF^{*}(\alpha,\beta)=G(\alpha,\beta)=0$, and \eqref{FFstar} still holds. Therefore, for complex vectors $h,k$,
\begin{gather}
|\langle FF^*h,k\rangle|
\le\langle FF^*|h|,|k|\rangle
\le\BoxC{w}{\mu}\langle G|h|,|k|\rangle
\le\BoxC{w}{\mu}\|G\|\|h\|_2\|k\|_2.
\label{eq:absolutecomparison}
\end{gather}
Taking supremum over all $h,k$ with $\|h\|_2=\|k\|_2=1$ proves the result. 
\end{proof}

We have the following more general family of matrices, which includes $F$ when $S=\varnothing$ and $F^*$ when $S=\{1,...,n\}$:

\begin{definition}
For a coordinate split $S$, define the oriented matrix $T_S$ by
\begin{gather*}
T_S((a,v),(u,b))=
\begin{cases}
\sqrt{w(a,v)w(u,b)}\,(\Hardy^*\mu)(a,b),&a\le u,\ b\le v,\\
0,&\text{otherwise}.
\end{cases}
\end{gather*}
\end{definition}

Note that $a,u\in\mathcal T_S$ and $b,v\in\mathcal T_{S^c}$, and the coordinates of each paired vertex are given according to \eqref{eq:block-assembly}. 

\begin{lemma}\label{lem:trace}
For a fixed split satisfying the corner condition \eqref{eq:corner}, every positive integer $m$ satisfies
\begin{gather*}
\Tr\left((T_ST_S^*)^m \right)\le K_S^m\Tr\left((F^*F)^m\right).
\end{gather*}
Consequently $\|T_S\|^2\le K_S\|F\|^2$. This result remains valid when the cumulative mass function $\Hardy^*\mu$ in $T_S$ and $F$ is replaced by any nonnegative function $g$ on $T^n$.
\end{lemma}

\begin{proof}
Let $g$ denote $\Hardy^*\mu$ or any nonnegative function replacing it in both matrices.

Fix $m\ge1$. All tuple indices are read modulo $m$: for each tuple $x$ we set $x_{i+m}=x_i$. 

For row labels $(a,v)$ and $(a',v')$, matrix multiplication gives
\begin{align*}
(T_ST_S^*)((a,v),(a',v'))
&=\sum_{u,b}T_S((a,v),(u,b))
\overline{T_S((a',v'),(u,b))}
\\
&=\sqrt{w(a,v)w(a',v')}
\sum_{\substack{a,a'\le u\\b\le v,v'}}
w(u,b)g(a,b)g(a',b),
\end{align*}
since the entries of $T$ are nonnegative. Writing $\alpha_i=(a_i,v_i)$ and $\beta_i=(u_i,b_i)$, we obtain
\begin{align*}
\Tr((T_ST_S^*)^m)
&=\sum_{\alpha_1,\ldots,\alpha_m}
\prod_{i=1}^m(T_ST_S^*)(\alpha_i,\alpha_{i+1})
\\
&=\sum_{\substack{\alpha_1,\ldots,\alpha_m\\
\beta_1,\ldots,\beta_m}}
\prod_{i=1}^m T_S(\alpha_i,\beta_i)T_S(\alpha_{i+1},\beta_i).
\end{align*}

By definition of $T_S$, the product above vanishes unless the tuple
$(a, b, u, v)$ belongs to the set $\Omega$ defined by
\begin{gather}
a_i,a_{i+1}\le u_i,\qquad
b_i\le v_i,v_{i+1}
\quad(1\le i\le m).
\label{eq:omega}
\end{gather}

For a tuple in $\Omega$, we rewrite the product above:
\begin{eqnarray*}
\prod_{i=1}^m T_S(\alpha_i,\beta_i)T_S(\alpha_{i+1},\beta_i) 
&=& 
\prod_{i=1}^m 
\sqrt{w(a_i,v_i)w(a_{i+1},v_{i+1})}\,
w(u_i,b_i)g(a_i,b_i)g(a_{i+1},b_i)\\
&=& 
\prod_{i=1}^m w(a_i,v_i)\,
w(u_i,b_i)g(a_i,b_i)g(a_{i+1},b_i)\\
&=:&
\Pi_0
\end{eqnarray*}
due to the cyclic convention.

Therefore,
\begin{gather*}
\Tr((T_ST_S^*)^m)=\sum_\Omega\Pi_0,
\end{gather*}

Next, we give two different representations for $\Tr((F^*F)^m)$. By the definition of $F$,
\begin{align*}
(F^*F)(\alpha,\alpha')
&=\sum_{\beta\ge\alpha,\alpha'}
F(\beta,\alpha)F(\beta,\alpha')
\\
&=\sqrt{w(\alpha)w(\alpha')}g(\alpha)g(\alpha')
\sum_{\beta\ge\alpha,\alpha'}w(\beta).
\end{align*}
Consequently,
\begin{align*}
\Tr((F^*F)^m)
&=\sum_{\alpha_1,\ldots,\alpha_m\in T^n}
\prod_{i=1}^m(F^*F)(\alpha_i,\alpha_{i+1})
\\
&=
\sum_{\alpha_1,\ldots,\alpha_m\in T^n}
\prod_{i=1}^m
\bigg(\sqrt{w(\alpha_i)w(\alpha_{i+1})}
g(\alpha_i)g(\alpha_{i+1})
\sum_{\beta_i\ge\alpha_i,\alpha_{i+1}}w(\beta_i)
\bigg)
\\
&=\sum_{\alpha_1,\ldots,\alpha_m\in T^n}
\prod_{i=1}^m
\bigg(
w(\alpha_i)g(\alpha_i)^2
\sum_{\beta_i\ge\alpha_i,\alpha_{i+1}}w(\beta_i)
\bigg),
\end{align*}
by the cyclic convention. For the first representation, choose
\begin{gather*}
\alpha_i=(a_i,b_i),\qquad
\beta_i=(u_i,v_{i+1}).
\end{gather*}
The conditions $\alpha_i,\alpha_{i+1}\le\beta_i$ become
\begin{gather*}
a_i,a_{i+1}\le u_i,\qquad
b_i,b_{i+1}\le v_{i+1},
\end{gather*}
which describe $\Omega$. Thus,
\begin{eqnarray*}
\Tr((F^*F)^m) 
&=& 
\sum_\Omega \prod_{i=1}^m
g(a_i,b_i)^2w(a_i,b_i)w(u_i,v_{i+1}) \  
=:\  \sum_\Omega\Pi_1,
\end{eqnarray*}

For the second representation, choose instead
\begin{gather*}
\alpha_i=(a_i,b_{i-1}),\qquad
\beta_i=(u_i,v_i),
\end{gather*}
Now the conditions become
\begin{gather*}
a_i,a_{i+1}\le u_i,\qquad
b_{i-1},b_i\le v_i,
\end{gather*}
again describing $\Omega$. Moreover, reindexing gives
\begin{gather*}
\prod_{i=1}^m g(a_i,b_{i-1})^2w(a_i,b_{i-1})
=\prod_{i=1}^m g(a_{i+1},b_i)^2w(a_{i+1},b_i).
\end{gather*}
Therefore, we also have
\begin{eqnarray*}
\Tr((F^*F)^m) 
&=& 
\sum_\Omega \prod_{i=1}^m
g(a_{i+1},b_i)^2w(a_{i+1},b_i)w(u_i,v_i)
\ =:\ \sum_\Omega\Pi_2,
\end{eqnarray*}

We now compare the three products. For every tuple in $\Omega$, the corner condition \eqref{eq:corner} gives
\begin{align*}
w(a_i,v_i)w(u_i,b_i)
&\le K_S\ w(a_i,b_i)w(u_i,v_i),
\\
w(a_{i+1},v_{i+1})w(u_i,b_i)
&\le K_S\ w(a_{i+1},b_i)w(u_i,v_{i+1}).
\end{align*}

Multiplying both inequalities over $i$ and using cyclicity on the left gives
\begin{align*}
\bigg(\prod_{i=1}^m w(a_i,v_i)w(u_i,b_i)\bigg)^2
&\le 
K_S^{2m}\prod_{i=1}^m
w(a_i,b_i)w(a_{i+1},b_i)
w(u_i,v_i)w(u_i,v_{i+1}).
\end{align*}

Multiplying by
$\prod_{i=1}^m g(a_i,b_i)^2g(a_{i+1},b_i)^2$
therefore yields
\begin{gather*}
\Pi_0^2\le (K_S^m\Pi_1)(K_S^m\Pi_2).
\end{gather*}

Summing over $\Omega$, and applying Cauchy--Schwarz now gives
\begin{align}\label{trace ineq}
\Tr((TT^*)^m)
=\sum_\Omega\Pi_0
&\le K_S^m\sum_\Omega\sqrt{\Pi_0\Pi_1}
\\
&\le K_S^m
\bigg(\sum_\Omega\Pi_0\bigg)^{1/2}
\bigg(\sum_\Omega\Pi_1\bigg)^{1/2} \notag
\\
&=K_S^m\Tr((F^*F)^m). \notag
\end{align}

Finally, recall that for a finite positive semidefinite matrix $M$ of size $N$, its eigenvalues $\lambda_1,\ldots,\lambda_N$ are nonnegative, and
\begin{gather*}
\|M\|^m=(\max_j \lambda_j )^m \le 
\sum_{j=1}^N\lambda_j^m=\Tr(M^m)
\le N\|M\|^m.
\end{gather*}
Taking $m$th roots and letting $m\to\infty$ gives
$\Tr(M^m)^{1/m}\to\|M\|$.

Since the matrices $T_ST_S^*$ and $F^*F$ are positive semidefinite, taking $m$th roots in \eqref{trace ineq} and passing to the limit gives
\begin{gather*}
\|T\|^2=\|TT^*\|
\le K\|F^*F\|=K\|F\|^2.
\end{gather*}
\end{proof}

\begin{proof}[Proof of Theorem~\ref{thm:corner}]
First, we prove that the matrices $T_S$ cover $G$. For an entry with $\alpha\sim\beta$, where $\alpha=(\alpha_1,\alpha_2,\dots), \beta=(\beta_1,\beta_2,\dots)$,  choose the coordinate split
\[
S=\{j:\alpha_j\le\beta_j\}.
\]
In the complement of $S$ one has $\beta_j<\alpha_j$. Thus $T_S(\alpha,\beta)=G(\alpha,\beta)$ for this choice of $S$. For $\alpha\nsim\beta$, $T_S(\alpha,\beta)=G(\alpha,\beta)=0$ for all the splits $S$. Thus, 

$$
0\le G(\alpha,\beta)\le\sum_{S\subseteq\{1,\ldots,n\}}T_S(\alpha,\beta), 
$$ 
which implies (with a similar argument as in \eqref{eq:absolutecomparison})
$$
\|G\|\le\sum_{S\subseteq\{1,\ldots,n\}}\|T_S\|.
$$ 
Now, lemmata ~\ref{lem:energy} and~\ref{lem:trace} imply
\[
\|G\|\le\sum_S\|T_S\|
\le\left(\sum_S\sqrt{K_S}\right)\|F\|
\le\left(\sum_S\sqrt{K_S}\right)
\sqrt{\BoxC{w}{\mu}\,\|G\|}.
\]
Using \eqref{eq:norms}, if $\|G\|=0$, \eqref{eq:directcorner}  holds trivially. Otherwise the display above gives \eqref{eq:directcorner}. 
For \eqref{eq:globalK}, the empty and full split $S$ have $K_S=1$, while any other split (in total  $2^n-2$) has constant $K_S\leq K$, so $\kappa(w)\le\bigl(2+(2^n-2)\sqrt K\bigr)^2\le 4^nK$
\end{proof}

\begin{proof}[Proof of Theorem~\ref{thm:product}] 
Note that product weights satisfy the unit corner condition \eqref{eq:unitcorner} as equality. Theorem~\ref{thm:corner} gives the upper bound $4^n$ for every
weight satisfying the unit corner condition \eqref{eq:unitcorner}. We prove sharpness by constructing a measure $\mu$ and weights $w$ so that for every $C<4^n$, $[w,\mu]_{\text{CE}}>C[w,\mu]_{\text{Box}}$. We use $w\equiv 1$, which are product weights and therefore satisfy the unit corner condition \eqref{eq:unitcorner}. To construct $\mu$, fix integers $1\le M<N$ and let $T$ be the chain
\[
\alpha_N\le\alpha_{N-1}\le\cdots\le\alpha_M.
\]

Define a measure on $T$ by
\[
\mu(\alpha_j)=
\begin{cases}
\dfrac1j-\dfrac1{j+1},&M\le j<N,\\[10pt]
\dfrac1N,&j=N.
\end{cases}
\]

Therefore, for some integer $M\leq k \leq N$,
$$
(\Hardy^*\mu)(\alpha_k)=\sum_{j=k}^N \mu(a_j)=\frac{1}{k}.
$$
Consequently,
\[
\BoxC{1}{\mu}
=
\sup_{M\le k\le N}\frac{\sum\limits_{a_j\le a_k}((\Hardy^*\mu)(\alpha_j))^2}{(\Hardy^*\mu)(\alpha_k)}
=
\max_{M\le k\le N}k\sum_{j=k}^N\frac1{j^2}
\le 1+\frac1M,
\]
where we used that
$\displaystyle\sum_{j=k}^N \frac{1}{j^2}\le \frac{1}{k^2}+\int_k^\infty \frac{dt}{t^2}=\frac{1}{k^2}+\frac{1}{k}$.

Define now
\[
f(\alpha_j)=
\begin{cases}
\dfrac{\sqrt{j(j+1)}}{\sqrt j+\sqrt{j+1}},&M\le j<N,\\[10pt]
\sqrt N,&j=N.
\end{cases}
\]
For $j<N$, we have
$$
f(\alpha_j)\mu(\alpha_j)=\frac{1}{\sqrt{j}}-\frac{1}{\sqrt{j+1}},
$$
while $f(\alpha_N)\mu(\alpha_N)=\frac{1}{\sqrt{N}}$. Thus, 
$$
(\Hardy^*(f\mu))(\alpha_j)=\sum_{i=j}^N f(\alpha_i)\mu(\alpha_i)=\frac{1}{\sqrt{j}}
$$
and hence
\[
\E_1[f\mu]=\sum_{j=M}^N\frac1j
\ \ \text{ and }\ \
\|f\|_{L^2(\mu)}^2
=1+\sum_{j=M}^{N-1}\frac1{(\sqrt j+\sqrt{j+1})^2}.
\]
The inequalities
\[
\frac{1}{4(j+1)}
\le \frac{1}{(\sqrt j+\sqrt{j+1})^2}
\le \frac{1}{4j}
\]
give
\[
1+\frac14\left(\E_1[f\mu]-\frac1M\right)
\le \|f\|_{L^2(\mu)}^2
\le 1+\frac14\left(\E_1[f\mu]-\frac1N\right).
\]
For fixed $M$, if $N\to\infty$, then $\E_1[f\mu]\to\infty$
and consequently 
\[
\frac{\E_1[f\mu]}{\|f\|_{L^2(\mu)}^2}
\longrightarrow 4.
\]

Combining this limit with the bound $\BoxC{1}{\mu}\le 1+\frac{1}{M}$, we have
\[
\liminf_{N\to\infty}
\frac{\E_1[f\mu]}{\BoxC{1}{\mu}\,\|f\|_{L^2(\mu)}^2}
\ge \frac{4}{1+\frac{1}{M}}.
\]

Now take the product of $n$ copies of this chain. Hence, we get an $n$-tree $T^n$ and we set for $\alpha=(\alpha_1, \dots, \alpha_n)$,
$$
\nu(\alpha)=\prod_{j=1}^{n}\mu(\alpha_j)
\ \ \text{ and }\ \ 
g(\alpha)=\prod_{j=1}^{n}f(\alpha_j),
$$
still with $w\equiv1$.

Through Fubini's theorem and the coordinate-wise  factoring of the function $g$ and the measure $\nu$ we get
\[
\BoxC{1}{\nu}=\BoxC{1}{\mu}^{\,n}.
\]

Likewise,
\[
\E_1[g\nu]=\bigl(\E_1[f\mu]\bigr)^n 
\ \ \text{ and }\ \ 
\|g\|_{L^2(\nu)}^2=\|f\|_{L^2(\mu)}^{2n}.
\]

Hence
\[
\frac{\CE{1}{\nu}}{\BoxC{1}{\nu}}
\ge
\frac{\E_1[g\nu]}{\BoxC{1}{\nu}\,\|g\|_{L^2(\nu)}^2}
=
\left(
\frac{\E_1[f\mu]}{\BoxC{1}{\mu}\,\|f\|_{L^2(\mu)}^2}
\right)^n
\ge
\left( \frac{4}{1+\frac{1}{M}} \right)^n.
\]

Given any $C<4^n$, first choose $M$ sufficiently large that
\[
\left(\frac{4}{1+\frac{1}{M}}\right)^n>C.
\]

Then, choose $N$ sufficiently large so to get \(\CE{1}{\nu}>C\,\BoxC{1}{\nu}\). Therefore, $4^n$ is sharp among all products weights and therefore for all weights satisfying the corner condition \eqref{eq:corner}.
\end{proof}

\begin{remark}\label{rem:bounded-perturbations}
Suppose $w_0$ is a product weight and
$a w_0\le w\le b w_0$, where $0<a\le b<\infty$. Then
\[
\CE{w}{\mu}\le4^n\frac{b}{a}\BoxC{w}{\mu}.
\]
Indeed, $\E_w[f\mu]\le b\E_{w_0}[f\mu]$ and
$\BoxC{w_0}{\mu}\le\frac{\BoxC{w}{\mu}}{a}$, so
Theorem~\ref{thm:product} applies. The same estimate holds when
$w_0$ satisfies the unit corner condition \eqref{eq:unitcorner},
by Theorem~\ref{thm:corner}. Thus the box characterization is
preserved under bounded multiplicative perturbations of product
weights and, more generally, of weights satisfying the unit
corner condition.
\end{remark}

\begin{remark} [Infinite vertex trees] The corner theorem extends to locally finite rooted $n$-trees with infinitely many vertices
(each vertex has finitely many children), for finite measures on the
full vertex set. The box constant uses the full descendant sums in
\eqref{eq:boxconstant}.

To pass to infinite depth, we cut the $n$-tree below depth $N$ in each coordinate and let $T^{n,N}$ be their full product.
Local finiteness ensures that $T^{n,N}$ is finite.
Let $\pi_N:T^n\to T^{n,N}$ act coordinate-wise, fixing surviving vertices and sending every deeper vertex to its depth-$N$ ancestor.
Set $\mu_N=(\pi_N)_*\mu$ for the corresponding push-forward measure:
\begin{gather*}
\mu_N(\eta)=\mu(\pi_N^{-1}(\{\eta\})).
\end{gather*}
For $f\in L^2(\mu)$, define
\begin{gather*}
f_N(\eta)=
\begin{cases}
\displaystyle
\frac{1}{\mu_N(\eta)}
\int_{\pi_N^{-1}(\{\eta\})}f\,d\mu,
&\mu_N(\eta)>0,\\
0,&\mu_N(\eta)=0.
\end{cases}
\end{gather*}
For every surviving vertex $\alpha$,
\begin{gather*}
\pi_N^{-1}\{\eta\in T^{n,N}:\eta\le\alpha\}
=\{\beta\in T^n:\beta\le\alpha\}.
\end{gather*}
This identity and Cauchy--Schwarz on each set $\pi_N^{-1}(\{\eta\})$ give
\begin{align}
(\Hardy^*\mu_N)(\alpha)
&=(\Hardy^*\mu)(\alpha),
\notag\\
(\Hardy^*(f_N\mu_N))(\alpha)
&=(\Hardy^*(f\mu))(\alpha),
\notag\\
\|f_N\|_{L^2(\mu_N)}^2
&\le\sum_{\eta\in T^{n,N}}
\int_{\pi_N^{-1}(\{\eta\})}|f|^2\,d\mu
=\|f\|_{L^2(\mu)}^2.
\label{eq:truncation-identities}
\end{align}

Put $w_N=w|_{T^{n,N}}$. For every $\gamma\in T^{n,N}$,
\begin{align*}
\sum_{\substack{\alpha\in T^{n,N}\\\alpha\le\gamma}}
w_N(\alpha)(\Hardy^*\mu_N)(\alpha)^2
\le\sum_{\alpha\le\gamma}
w(\alpha)(\Hardy^*\mu)(\alpha)^2
\le\BoxC{w}{\mu}\,(\Hardy^*\mu_N)(\gamma).
\end{align*}

Hence
\begin{gather*}
\BoxC{w_N}{\mu_N}\le\BoxC{w}{\mu}.
\end{gather*}

The finite estimates give
\begin{align*}
\E_{w_N}[f_N\mu_N]
\le 
\kappa(w_N)\BoxC{w_N}{\mu_N}\|f_N\|_{L^2(\mu_N)}^2
\le 
\kappa(w)\BoxC{w}{\mu}\|f\|_{L^2(\mu)}^2.
\end{align*}
By \eqref{eq:truncation-identities},
\begin{gather*}
\E_{w_N}[f_N\mu_N]
=\sum_{\alpha\in T^{n,N}}
w(\alpha)|(\Hardy^*(f\mu))(\alpha)|^2
\end{gather*} which increases to $\E_w[f\mu]$ by Monotone convergence. Therefore, 
\begin{gather*}
\BoxC{w}{\mu}\le\CE{w}{\mu}
\le \kappa(w)\BoxC{w}{\mu}.
\end{gather*}
\end{remark}

\section{Capacity and Carleson measures for Dirichlet spaces on the polydisc}
\label{sec:capacity}
We now apply the unweighted embedding theorem to the classical
Dirichlet space on the polydisc. The analytic-to-dyadic reduction in
\cite[Theorem~2.2]{MPVZ} applies in every
dimension. Our theorem supplies the box characterization of the
discrete embeddings arising in that reduction.

Define the Whitney tiles by
\begin{gather*}
q_\theta(I)=\left\{re^{2\pi it}:t\in I,\ 1-|I|\le r<1-\frac{|I|}{2}\right\}.
\end{gather*}
The origin belongs only to the root tile and the root tent. With the
assigned edges, the tiles form a Borel partition of $\mathbb D$, and \eqref{Stheta} gives
\begin{gather*}
S_\theta(I)=\bigsqcup_{\substack{J\in\mathcal D^\theta\\J\subseteq I}}
q_\theta(J),
\end{gather*}

Regard $\mathcal R^\omega$ as an infinite binary $n$-tree $T_\omega^n$.
For $\alpha=I_1\times\cdots\times I_n\in\mathcal R^\omega$, set
$q_\omega(\alpha)=\prod\limits_{j=1}^n q_{\theta_j}(I_j)$ and define the
discrete Whitney-tile measure by
\begin{gather*}
\sigma_\omega(\alpha)=\nu(q_\omega(\alpha)).
\end{gather*}

Thus $\sigma_\omega$ has total mass $\nu(\mathbb D^n)$; its mass at
one vertex is the mass of a Whitney tile, and its cumulative mass is
the mass of the corresponding tent:
$$
(\Hardy^*\sigma_\omega)(\alpha)
=\sum_{\beta\le\alpha}\nu(q_\omega(\beta))
=\nu(S_\omega(\alpha)).
$$

In particular, the box constant
in \eqref{eq:polydisc-uniform-box-test} satisfies
\begin{gather}
[\nu]_{\mathrm{Box}}=\sup_\omega\BoxC{1}{\sigma_\omega}.
\label{eq:polydisc-tile-tent-masses}
\end{gather}

We now state the comparison between the analytic and dyadic embedding constants. Let
the infinite vertex embedding
\begin{gather*}
\sum_{\alpha\in T_\omega^n}
|\Hardy^*(g\sigma_\omega)(\alpha)|^2
\le \CE{1}{\sigma_\omega}
\sum_{\beta\in T_\omega^n}|g(\beta)|^2\sigma_\omega(\beta).
\end{gather*}

If we specialize \cite[Theorem~2.2]{MPVZ} to
$\vec s=(1,\ldots,1)$ gives constants $a_n,b_n\in(0,\infty)$,
depending only on $n$, such
that
\begin{gather}
[\nu]_{\mathrm{CM}}
\le a_n\sup_{\omega\in\mathbb T^n}\CE{1}{\sigma_\omega}\
\ \text{and}\ \ 
\sup_{\omega\in\mathbb T^n}\CE{1}{\sigma_\omega}
\le b_n[\nu]_{\mathrm{CM}}.
\label{eq:polydisc-published-reduction}
\end{gather}

These are the constants in
Theorem~\ref{thm:polydisc-dirichlet}. They absorb the equivalence
between \eqref{eq:polydisc-dirichlet-norm} and the norm associated with
the logarithmic reproducing kernel in \cite[Section~2.2]{MPVZ}. This reduction is independent of the box characterization on $n$-trees. 

\begin{proof}[Proof of Theorem~\ref{thm:polydisc-dirichlet}]
Fix a rotation $\omega$. Since $w\equiv1$ is a product weight, by
Theorem~\ref{thm:product}, extended to the infinite vertex trees described in Remark \ref{eq:boxconstant},
\begin{gather*}
\CE{1}{\sigma_\omega}\le4^n\BoxC{1}{\sigma_\omega}.
\end{gather*}

Testing the embedding with the characteristic function of a descendant set, we get that $\BoxC{1}{\sigma_\omega}\le \CE{1}{\sigma_\omega}$. Taking supremum over all rotations and using \eqref{eq:polydisc-tile-tent-masses},
\begin{gather*}
[\nu]_{\mathrm{Box}}
\le\sup_\omega \CE{1}{\sigma_\omega}
\le4^n[\nu]_{\mathrm{Box}}.
\end{gather*}

Combining these inequalities with
\eqref{eq:polydisc-published-reduction}, we get the bounds in \eqref{eq:polydisc-main-bounds}.
\end{proof}

\begin{comment}
\begin{remark}[Product Dirichlet weights]\label{rem:other-dirichlet-weights}
The same argument applies to the weighted Dirichlet spaces
$H_{\vec s}(\mathbb D^n)$ considered in
\cite[Theorem~2.2]{MPVZ}, provided
$1-\varepsilon_n\le s_j\le1$ for every $j$, with
$\varepsilon_n>0$ as in that theorem.
Indeed, the analytic-to-dyadic reduction produces the product weight
\[
w_{\vec s}(R)=\prod_{j=1}^n |I_j|^{s_j-1},
\qquad R=I_1\times\cdots\times I_n.
\]
Our product-weight theorem therefore shows that
$H_{\vec s}(\mathbb D^n)\hookrightarrow L^2(\nu)$
if and only if there is a finite constant $C$ such that
\[
\sum_{\substack{R\in\mathcal R^\omega\\R\subseteq Q}}
w_{\vec s}(R)\,\nu(S_\omega(R))^2
\le C\,\nu(S_\omega(Q))
\]
for every rotation $\omega$ and every
$Q\in\mathcal R^\omega$.
The embedding and testing constants are comparable, with
comparison constants depending only on $n$ and $\vec s$.
\end{remark}
\end{comment}

\begin{proof}[Proof of Corollary~\ref{cor:polydisc-composition}]
The definition of $\nu_\Phi$ gives
\begin{gather*}
\|C_\Phi f\|_{A^2_\beta(\mathbb D^n)}^2
=\int_{\mathbb D^n}|f(\Phi(z))|^2\,dV_\beta(z)
=\int_{\mathbb D^n}|f|^2\,d\nu_\Phi.
\end{gather*}
Thus $\|C_\Phi\|^2=[\nu_\Phi]_{\mathrm{CM}}$, and
Theorem~\ref{thm:polydisc-dirichlet} now gives the result we want to prove.
\end{proof}

\subsection*{Capacity and Carleson Embedding}
We relate the embedding inequalities to discrete capacity. The capacity comparison itself holds for a general family of measurable
sets and does not require product geometry.

Let $X$ be a measurable space, let $\mu$ be a finite nonnegative measure
on $X$, and let $\mathcal J$ be a finite or countable index set. For each
$\alpha\in\mathcal J$, choose a measurable set $D_\alpha\subseteq X$,
and let $w:\mathcal J\to[0,\infty)$. For nonnegative
$h\in\ell^2(\mathcal J)$, define
\begin{gather*}
A_wh(z)=\sum_{\substack{\alpha\in\mathcal J\\ z\in D_\alpha}}
\sqrt{w(\alpha)}h(\alpha)=\sum_{\alpha\in\mathcal J}
\sqrt{w(\alpha)}h(\alpha)\mathbf{1}_{D_\alpha}.
\end{gather*}
The sum may be infinite. For a measurable set $E\subseteq X$, define
\begin{gather}
\operatorname{Cap}_w(E)
=\inf\bigg\{\sum_{\alpha\in\mathcal J}h(\alpha)^2:
h\in\ell^2(\mathcal J),\ h\ge0,\ A_wh\ge1\text{ on }E\bigg\}.
\label{eq:weighted-capacity}
\end{gather}

We call the functions $h\in\ell^2$ admissible for $E$ if $h\ge 0$, $A_wh\ge 1$ on $E$. The above infimum is $+\infty$, if the admissible class is empty, and
$\operatorname{Cap}_w(\varnothing)=0$. This capacity depends on $w$ and
the sets $D_\alpha$, independently of $\mu$. Set
\begin{gather*}
C_{\mathrm{cap}}(w,\mu)
=\sup_{E\subseteq X\text{ measurable}}
\frac{\mu(E)}{\operatorname{Cap}_w(E)},
\end{gather*}
with the conventions $0/0=0$, $a/0=+\infty$ for $a>0$, and
$a/\infty=0$ for $a\neq 0$. Equivalently,
$C_{\mathrm{cap}}(w,\mu)=\inf\{C>0:\mu(E)\le C\operatorname{Cap}_w(E)\}$ for every measurable $E$.

We define the \textit{hereditary} Carleson constant as the least
constant in
\begin{gather}
\sum_{\alpha\in\mathcal J}w(\alpha)\mu(E\cap D_\alpha)^2
\le [w,\mu]_{\mathrm{HC}}\mu(E)
\qquad(E\subseteq X\text{ measurable}).
\label{eq:capacity-restricted-energy}
\end{gather}

\begin{lemma}\label{lem:capacity-restricted}
For every such \(\{D_\alpha\}_{\alpha\in\mathcal J}\) and every nonnegative weight,
\begin{gather}
C_{\mathrm{cap}}(w,\mu)
\le [w,\mu]_{\mathrm{HC}}
\le4C_{\mathrm{cap}}(w,\mu).
\label{eq:capacity-hc-comparison}
\end{gather}
\end{lemma}

\begin{proof}
Suppose that $[w,\mu]_{\mathrm{HC}}<\infty$. If $h$ is
admissible for $E$ in \eqref{eq:weighted-capacity}, Tonelli's theorem
and Cauchy--Schwarz give
\begin{align*}
\mu(E)&\le\int_E A_wh\,d\mu
=\sum_\alpha\sqrt{w(\alpha)}h(\alpha)\mu(E\cap D_\alpha)\\
&\le\|h\|_2\bigg(\sum_\alpha
w(\alpha)\mu(E\cap D_\alpha)^2\bigg)^{1/2}
\le\|h\|_2\sqrt{[w,\mu]_{\mathrm{HC}}\mu(E)}.
\end{align*}
For $\mu(E)>0$, the above yields
$\mu(E)\le [w,\mu]_{\mathrm{HC}}\|h\|_2^2$. Taking the infimum over admissible $h$ proves that
$\frac{\mu(E)}{\operatorname{Cap}_w(E)}\le [w,\mu]_{\mathrm{HC}}$ and then taking supremum over all $E$, we have the first inequality in \eqref{eq:capacity-hc-comparison}. Sets of zero
measure and sets of infinite capacity contribute zero to the supremum.
If $[w,\mu]_{\mathrm{HC}}=+\infty$, this inequality is automatic.

For the reverse inequality, take $C>0$ such that
$\mu(F)\le C\operatorname{Cap}_w(F)$ for every measurable $F$.
Let $g\ge0$ have finite support in $\mathcal J$. For $t>0$, the function
$\dfrac{g}{t}$ is admissible for 
$$
F_t=\{z:A_wg(z)>t\},
$$ 
so
$$
\mu(F_t)\le C \operatorname{Cap}_w(F_t)\le C\frac{\|g\|_2^2}{t^2}.
$$

For $g\ne0$, we integrate over a measurable $E$ with $\mu(E)>0$ and we obtain
\begin{align*}
\int_E A_wg\,d\mu
&=
\int_E \int_0^\infty \mathbf{1}_{F_t}(t)dtd\mu(z)
=
\int_0^\infty \mu(E\cap F_t) dt \le
\int_0^\infty
\min\bigg\{\mu(E),\frac{C\|g\|_2^2}{t^2}\bigg\}\,dt \notag
\\
&=
\int_0^{t_0}
\min\bigg\{\mu(E),\frac{C\|g\|_2^2}{t^2}\bigg\}\,dt
+
\int_{t_0}^\infty
\min\bigg\{\mu(E),\frac{C\|g\|_2^2}{t^2}\bigg\}\,dt
=
2\sqrt{C\mu(E)}\|g\|_2,
\label{eq:capacity-layercake}
\end{align*}
where $t_0=\sqrt{\frac{C}{\mu(E)}}\|g\|_2$. If $g=0$ or $\mu(E)=0$, then the result is
immediate.

Fix a finite set $\mathcal F\subseteq\mathcal J$ and set
$g(\alpha)=\sqrt{w(\alpha)}\mu(E\cap D_\alpha)$ on $\mathcal F$, with
$g=0$ elsewhere. Set
$L=\sum\limits_{\alpha\in\mathcal F}w(\alpha)\mu(E\cap D_\alpha)^2$. Then $\int_E A_wg\,d\mu=L$ and
$\|g\|_2=\sqrt L$. By the estimate above, $L\le2\sqrt{C\mu(E)L}$. This implies $L\le4C\mu(E)$, when $L\neq 0$. For $L=0$, this inequality is trivial. Taking the supremum over finite
$\mathcal F$, we get $[w,\mu]_{\mathrm{HC}}
\le4C$. Letting $C$ decrease to $C_{\mathrm{cap}}(w,\mu)$ proves $[w,\mu]_{\mathrm{HC}}
\le4C_{\mathrm{cap}}(w,\mu)$. If this constant is infinite, the assertion is
automatic.
\end{proof}

On an $n$-tree, take $D_\alpha=\{\beta\in T^n:\beta\le\alpha\}$. Then
$A_wh=\Hardy(\sqrt w\,h)$ and 
$$
\E_w[\mu|_E]=\sum\limits_{\alpha\in T^n }w(\alpha)|\sum\limits_{D_\alpha}\mu(\beta)\mathbf{1}_E(\beta)|^2 = \sum\limits_{\alpha\in T^n}w(\alpha)\mu(E\cap D_\alpha)^2.
$$
Thus,
\eqref{eq:capacity-restricted-energy} becomes
\begin{gather*}
\E_w[\mu|_E]\le\HC{w}{\mu}\,\mu(E)
\qquad(E\subseteq T^n).
\end{gather*}
For $w\equiv1$, \eqref{eq:weighted-capacity} is the usual discrete
Hardy capacity $\inf\{\|h\|_2^2:h\ge0,\ \Hardy h\ge1\text{ on }E\}$.
The box characterization therefore gives the following capacitary
criterion.

\begin{theorem}\label{thm:capacity-characterization}
Suppose that $w$ satisfies the corner hypotheses of
Theorem~\ref{thm:corner}. Then
\begin{gather*}
C_{\mathrm{cap}}(w,\mu)\le\CE{w}{\mu}
\le4\kappa(w) C_{\mathrm{cap}}(w,\mu).
\end{gather*}
In particular, if $w$ satisfies the unit corner condition
\eqref{eq:unitcorner}, then
\begin{gather*}
\CE{w}{\mu}\le4^{n+1}C_{\mathrm{cap}}(w,\mu).
\end{gather*}
\end{theorem}

\begin{proof}
Testing the embedding with $\ind_E$ gives
$\HC{w}{\mu}\le\CE{w}{\mu}$. Also, for a vertex $\gamma$,
$$
\E_w[\mu|_{D_\gamma}] = \sum_{\alpha}w(\alpha)\mu(D_\gamma\cap D_\alpha)^2 
\ge
\sum_{\alpha\le\gamma}w(\alpha)\mu(D_\alpha)^2
=
\sum_{\alpha\le\gamma}w(\alpha)((\Hardy^*\mu)(\alpha))^2
$$

Therefore, $\BoxC{w}{\mu}\le\HC{w}{\mu}$. Lemma~\ref{lem:capacity-restricted},
Theorem~\ref{thm:corner}, now give
\begin{align*}
C_{\mathrm{cap}}(w,\mu)
\le
\HC{w}{\mu}\le\CE{w}{\mu}
\le 
\kappa(w)\BoxC{w}{\mu}
\le 
\kappa(w)\HC{w}{\mu}
\le
4\kappa(w) C_{\mathrm{cap}}(w,\mu).
\end{align*}
If $w$ satisfies the unit corner condition, then $\kappa(w)=4^n$.
\end{proof}

\begin{remark}\label{rem:infinite-capacity}
The results extend to locally finite infinite trees, for finite measures $\mu$, with the same constants. Truncate each coordinate, restrict $w$, and push $\mu$ onto the truncation. This preserves cumulative masses at survived vertices, does not increase the box constant, and keeps the weight hypotheses with uniform coefficient $\kappa(w)$. For finitely supported $h$, the Hardy integral is unchanged by sufficiently deep truncations, so density and duality give the full embedding. Since the vertex set is countable, Lemma~\ref{lem:capacity-restricted} also gives the corresponding capacity comparisons.
\end{remark}

Theorem \ref{thm:capacity-characterization} can be used to prove the following: 

\begin{corollary}\label{cor:polydisc-capacity}
For each rotation $\omega$, let $\operatorname{Cap}_\omega$ be the
capacity on the vertex set $T_\omega^n$,
defined by \eqref{eq:weighted-capacity} with
$D_\alpha=\{\beta:\beta\le\alpha\}$. Set
\begin{gather*}
[\nu]_{\mathrm{cap,dyad}}
=\sup_{\omega\in\mathbb T^n}C_{\mathrm{cap}}(1,\sigma_\omega).
\end{gather*}
Then $\nu$ is a Carleson measure for $\mathcal D(\mathbb D^n)$ if and
only if $[\nu]_{\mathrm{cap,dyad}}<\infty$. More precisely,
\begin{gather*}
\frac{1}{b_n}[\nu]_{\mathrm{cap,dyad}}
\le[\nu]_{\mathrm{CM}}
\le a_n4^{n+1}[\nu]_{\mathrm{cap,dyad}}.
\end{gather*}
\end{corollary}

\begin{proof}
The infinite vertex case of
Theorem~\ref{thm:capacity-characterization} gives, for each rotation,
\begin{gather*}
C_{\mathrm{cap}}(1,\sigma_\omega)
\le \CE{1}{\sigma_\omega}
\le4^{n+1}C_{\mathrm{cap}}(1,\sigma_\omega).
\end{gather*}
Take the supremum over $\omega$ and apply
\eqref{eq:polydisc-published-reduction}.
\end{proof}

\subsection*{Logarithmic Dirichlet spaces}

We can use the corner theorem to a modification of the Dirichlet
space.
Fix $0<\tau<1$. For $f(z)=\sum_{k\in\mathbb N_0^n}a_kz^k$, set
\begin{gather*}
\|f\|_{\mathcal D_\tau^{\log}}^2
=\sum_{k\in\mathbb N_0^n}|a_k|^2
\prod_{j=1}^n(k_j+1)
\biggl(1+\sum_{j=1}^n\log(k_j+1)\biggr)^\tau.
\end{gather*}
Let $\mathcal D_\tau^{\log}(\mathbb D^n)$ consist of the holomorphic
functions such that this sum is finite.
Note that for $\tau=0$, the formula gives the classical Dirichlet norm.

We use the grids, tents, and tile measures introduced above.
If $d_j(\alpha)$ is the depth of the $j$th coordinate of $\alpha$,
with depth zero at the root, put
\begin{gather*}
W_\tau(\alpha)=
\biggl(1+\sum_{j=1}^n d_j(\alpha)\biggr)^{-\tau}.
\end{gather*}

\begin{theorem}[Dirichlet spaces with a joint logarithmic weight]
\label{thm:logarithmic-dirichlet}
Let $n\ge1$, $0<\tau<1$, and let $\nu$ be a finite nonnegative
Borel measure on $\mathbb D^n$. Then
$$
\int_{\mathbb D^n}|f|^2\,d\nu
\le C\|f\|_{\mathcal D_\tau^{\log}}^2
\qquad(f\in\mathcal D_\tau^{\log}(\mathbb D^n))
$$
holds with a finite constant if and only if
\begin{gather*}
B:=
\sup_{\omega\in\mathbb T^n}
\sup_{\substack{\gamma\in T_\omega^n\\
                \nu(S_\omega(\gamma))>0}}
\frac{\displaystyle\sum_{\alpha\le\gamma}
       W_\tau(\alpha)\nu(S_\omega(\alpha))^2}
     {\nu(S_\omega(\gamma))}
<\infty.
\end{gather*}
The least embedding constant is comparable to $B$, with constants
that depend only on $n$ and $\tau$.
Equivalently,
\begin{gather*}
\sup_{\omega\in\mathbb T^n}
C_{\mathrm{cap}}(W_\tau,\sigma_\omega)<\infty,
\end{gather*}
where $\sigma_\omega(\alpha)=\nu(q_\omega(\alpha))$.
\end{theorem}

For $0\leq t\leq1$, define
\[
L_t(q)=\sum_{k=0}^\infty
\frac{q^k}{(k+1)^{1+t}},
\qquad q\in\mathbb D.
\]

We first establish the following lemma, whose proof is left to the reader. 

\begin{lemma}\label{3.6}
Let
\(
\delta=\min(1,|1-q|)
\text{ and }
D=\left\lfloor\log_2\frac1\delta\right\rfloor.
\)
Then, uniformly for $0\leq t\leq1$,
\begin{equation}\label{1}
\operatorname{Re}L_t(q)
\approx
\operatorname{Re}L_t(1-\delta)
\approx
\sum_{m=0}^{D}2^{-tm}.
\end{equation}
Furthermore,
\begin{equation}\label{2}
\operatorname{Re}L_t(q)\geq c>0
\ \ \text{ and }\ \
|\operatorname{Im}L_t(q)|\leq C,
\end{equation}
where the constants are independent of $q$ and $t$.
\end{lemma}

We now define
\begin{equation}\label{5}
K_\tau(z,\zeta)
=
\tau\int_0^1t^{\tau-1}
\prod_{j=1}^n
\bigl(A_n+L_t(z_j\overline{\zeta_j})\bigr)\,dt,
\end{equation}
where $A_n>0$ sufficiently big, so that
\[
\left|
\arg\bigl(A_n+L_t(q)\bigr)
\right|
\leq \frac{\pi}{4n}.
\]
The argument of the product in \eqref{5} consequently has absolute
value at most $\pi/4$. It follows that
\begin{equation}\label{6}
|K_\tau(z,\zeta)|
\lesssim_n
\operatorname{Re}K_\tau(z,\zeta).
\end{equation}

We now find the coefficients of this kernel. For $k_j\geq1$,
the coefficient of $q^{k_j}$ in $A_n+L_t(q)$ is
$(k_j+1)^{-1-t}$. For $k_j=0$, the coefficient is $A_n+1>0$. Hence, the coefficient $c_k$ of $z^k\overline{\zeta}^k$ in
$K_\tau$ satisfies
\[
c_k
\approx_n
\frac1{\prod_{j=1}^n(k_j+1)}
\,
\tau\int_0^1t^{\tau-1}e^{-t\sum_{j=1}^n\log(k_j+1)}\,dt
\approx_{n,\tau}
\frac1{
\displaystyle
\prod_{j=1}^n(k_j+1)
\Big(1+\sum_{j=1}^n\log(k_j+1)\Big)^\tau
},
\]
due to the elementary estimate
\begin{equation}\label{7}
\tau\int_0^1t^{\tau-1}e^{-tM}\,dt
\approx_\tau
(1+M)^{-\tau},
\text{ for } M\geq0.
\end{equation}

The kernel $K_\tau$ thus defines a reproducing-kernel Hilbert
space whose norm is equivalent to \(\|\cdot\|_{\mathcal D_\tau^{\log}}\). Notice that the constants $A_n$ only
affect coefficients for which at least one $k_j=0$, and these
changes are by fixed multiplicative constants.

Consequently, $\|\cdot\|_{\mathcal D_\tau^{\log}}^2\approx\Vert T_{K_\tau}\Vert$, where $\Vert T_{K_\tau}\Vert$ is the norm of the operator:
\[
T_{K_\tau}g(z)
=
\int_{\mathbb D^n}K_\tau(z,\zeta)g(\zeta)\,d\nu(\zeta).
\]

We initially assume that $\nu$ is supported in a compact
subpolydisc, so that the kernel is bounded.

\begin{lemma}\label{3.7}
For each rotation $\omega\in\mathbb T^n$, define
\[
P_\omega(z,\zeta)
=
\sum_{\alpha\in T_\omega^n}
W_\tau(\alpha)
\mathbf 1_{S_\omega(\alpha)}(z)
\mathbf 1_{S_\omega(\alpha)}(\zeta).
\]
Then, for every $\omega$
\begin{equation}\label{11}
P_\omega(z,\zeta)
\lesssim_{n,\tau}
\operatorname{Re}K_\tau(z,\zeta)
\lesssim_{n,\tau}
\int_{\mathbb T^n}P_\omega(z,\zeta)\,d\omega.
\end{equation}
\end{lemma}

\begin{proof} First observe that
\[
\widetilde W_\tau(\alpha)
=
\tau\int_0^1
t^{\tau-1}
2^{-t\sum_{j=1}^nd_j(\alpha)}\,dt
\]
satisfies
\[
\widetilde W_\tau(\alpha)
\approx_\tau
\left(1+\sum_{j=1}^nd_j(\alpha)\right)^{-\tau}
=
W_\tau(\alpha).
\] by \eqref{7}, with $M=\left(\sum\limits_{j=1}^nd_j(\alpha)\right)\log2$. Then, for each coordinate, set 
\[
\delta_j
=
\min\left(1,\left|1-z_j\overline{\zeta_j}\right|\right)
\text{ and }
D_j=
\left\lfloor\log_2\frac1{\delta_j}\right\rfloor.
\]

Lemma \ref{3.6} gives
\begin{equation*}\label{14}
\operatorname{Re}L_t(z_j\overline{\zeta_j})
\approx
\sum_{m=0}^{D_j}2^{-tm},
\end{equation*}
which implies

\begin{equation}\label{16}
\operatorname{Re}K_\tau(z,\zeta)
\approx_{n,\tau}
\tau\int_0^1t^{\tau-1}
\prod_{j=1}^n
\left(\sum_{m=0}^{D_j}2^{-tm}\right)\,dt
\approx_{n,\tau}
\sum_{\substack{0\leq m_j\leq D_j\\1\leq j\leq n}}
\left(1+\sum_{j=1}^nm_j\right)^{-\tau}.
\end{equation}

Now, fix a rotation $\omega$. If $z$ and $\zeta$ belong to the same
$S_\omega(\alpha)$, then in each coordinate their boundary
directions lie in a common dyadic interval at the scale
determined by $d_j(\alpha)$. This is possible only when that
interval is sufficiently large relative to their separation.
In particular,
\[
d_j(\alpha)\leq D_j+C
\]
for a fixed absolute constant $C$.

At any fixed coordinate depth, there is at most a bounded
number of dyadic intervals that can contain both points. Thus
\[
P_\omega(z,\zeta) \lesssim_{n,\tau} \sum_{\substack{0\leq m_j\leq D_j+C\\1\leq j\leq n}}
\left(1+\sum_{j=1}^nm_j\right)^{-\tau}
\lesssim_{n,\tau}
\operatorname{Re}K_\tau(z,\zeta).
\]
by \eqref{16}, as adding $C$ extra levels it only changes this sum by a constant factor. This proves the left-hand side inequality of \eqref{11}.

For the other inequality, fix $j$ and $m$ so that
\(
2^{-m}\geq C\delta_j.
\)
If we rotate the dyadic grid uniformly, the probability that
the two boundary directions corresponding to $z_j$ and
$\zeta_j$ lie in the same dyadic interval at this scale is
bounded below by a positive constant. This happens since the two
directions are separated by at most a fixed fraction of the
interval length, and a uniformly rotated dyadic grid has a
uniformly positive probability of placing neither point across
a boundary.

The rotations in the $n$ coordinates are independent. Hence,
for every multi-index $m=(m_1,\ldots,m_n)$ satisfying
\[
0\leq m_j\leq D_j-C,
\]
the probability that the corresponding dyadic intervals are
common in all coordinates is bounded below by a constant
depending only on $n$.

It follows that
\[
\int_{\mathbb T^n}P_\omega(z,\zeta)\,d\omega
\gtrsim_{n,\tau}
\sum_{\substack{0\leq m_j\leq D_j-C\\1\leq j\leq n}}
\left(1+\sum_{j=1}^nm_j\right)^{-\tau}
\gtrsim_{n,\tau}
\operatorname{Re}K_\tau(z,\zeta)
,
\]
where the last inequality follows from \eqref{16}. This proves the right-hand side of \eqref{11}. 
\end{proof}

\begin{proof}[Proof of Theorem \ref{thm:logarithmic-dirichlet}] 
For each $\omega$, let $T_{P_\omega}$ be the integral operator
on $L^2(\nu)$ with kernel $P_\omega$:
\[
T_{P_\omega}g(z)
=
\int_{\mathbb D^n}P_\omega(z,\zeta)g(\zeta)\,d\nu(\zeta).
\]

Since $P_\omega$ is nonnegative and symmetric, its operator
norm can be computed by taking the supremum over nonnegative
functions $g$. The same is true for the positive kernel
$\operatorname{Re}K_\tau$. Moreover, by \eqref{6}, the operator norm
associated with $K_\tau$ is comparable to the one associated
with $\operatorname{Re}K_\tau$. Consequently, \eqref{11} gives
\[
\|T_{P_\omega}\|
\lesssim_{n,\tau}
\|T_{K_\tau}\|
\]
for every $\omega$, and
\[
\|T_{K_\tau}\|
\lesssim_{n,\tau}
\left\|
\int_{\mathbb T^n}T_{P_\omega}\,d\omega
\right\|
\leq
\sup_{\omega\in\mathbb T^n}\|T_{P_\omega}\|.
\]
Thus
\begin{equation}\label{20}
\|T_{K_\tau}\|
\approx_{n,\tau}
\sup_{\omega\in\mathbb T^n}\|T_{P_\omega}\|.
\end{equation}

Fix $\omega$. 
Set $q_\alpha=q_\omega(\alpha)$ and recall that
\[
S_\omega(\alpha)=\bigcup_{\beta\leq\alpha}q_\beta
\ \ \text{ and }\ \
\sigma_\omega(\beta)=\nu(q_\beta).
\]
The sets $q_\beta$ form a disjoint partition, up to sets of
$\nu$-measure zero. Hence, for any $g\in L^2(\nu)$,
\begin{equation}\label{21}
\int_{S_\omega(\alpha)}g\,d\nu
=
\sum_{\beta\leq\alpha}
\int_{q_\beta}g\,d\nu.
\end{equation}

Using the definition of $P_\omega$, 
\begin{align}\label{22}
\langle T_{P_\omega}g,g\rangle \notag
=
\int_{\mathbb D^n}\int_{\mathbb D^n}
P_\omega(z,\zeta)
g(\zeta)\overline{g(z)}
\,d\nu(\zeta)d\nu(z)
&=
\sum_{\alpha\in T_\omega^n}
W_\tau(\alpha)
\left|
\int_{S_\alpha}g\,d\nu
\right|^2\\
&=
\sum_{\alpha\in T_\omega^n}
W_\tau(\alpha)
\left|
\sum_{\beta\leq\alpha}
\int_{q_\beta}g\,d\nu
\right|^2.
\end{align}

Define
\(
f(\beta)=
\displaystyle
\frac1{\sigma(\beta)}
\int_{q_\beta}g\,d\nu,
\) if \(\sigma_\omega(\beta)>0\), and \(0\) otherwise. Then
\[
\int_{q_\beta}g\,d\nu=f(\beta)\sigma_\omega(\beta),
\]
and Jensen's inequality gives
\[
\begin{aligned}
\sum_\beta |f(\beta)|^2\sigma_\omega(\beta)
&=
\sum_{\beta:\sigma_\omega(\beta)>0}
\frac1{\sigma_\omega(\beta)}
\left|\int_{q_\beta}g\,d\nu\right|^2
\leq
\sum_\beta\int_{q_\beta}|g|^2\,d\nu
=
\int_{\mathbb D^n}|g|^2\,d\nu.
\end{aligned}
\]

Substituting these into \eqref{22}, we get
\[
\langle T_{P_\omega}g,g\rangle
\leq
[W_\tau,\sigma_\omega]_{\rm CE}
\|g\|_{L^2(\nu)}^2,
\]
where
\[
[W_\tau,\sigma_\omega]_{\rm CE}
=
\sup_{f\ne0}
\frac{
\displaystyle
\sum_\alpha W_\tau(\alpha)
\left|
\sum_{\beta\leq\alpha}f(\beta)\sigma_\omega(\beta)
\right|^2
}{
\displaystyle
\sum_\beta |f(\beta)|^2\sigma_\omega(\beta)
}.
\]
This proves
\(
\|T_{P_\omega}\|
\leq
[W_\tau,\sigma_\omega]_{\rm CE}.
\)

For the reverse inequality, take an  finitely
supported function $f$ on the tree and define a function $g$
that is constant on each tile \(
g(z)=f(\beta), \text{ if }z\in q_\beta.
\)
Then
\[
\|g\|_{L^2(\nu)}^2
=
\sum_\beta |f(\beta)|^2\sigma_\omega(\beta),
\]
and, by \eqref{21},
\[
\int_{S_\omega(\alpha)}g\,d\nu
=
\sum_{\beta\leq\alpha}f(\beta)\sigma_\omega(\beta).
\]
Therefore,
\[
\langle T_{P_\omega}g,g\rangle
=
\sum_\alpha W_\tau(\alpha)
\left|
\sum_{\beta\leq\alpha}f(\beta)\sigma_\omega(\beta)
\right|^2.
\]

Taking the supremum over finitely supported $f$, then
we get the reverse inequality. Thus
\(
\|T_{P_\omega}\|
=
[W_\tau,\sigma_\omega]_{\rm CE}.
\) 
Combining it with \eqref{20}, we conclude that
\begin{equation}\label{27}
\|T_{K_\tau}\|
\approx_{n,\tau}
\sup_{\omega\in\mathbb T^n}
[W_\tau,\sigma_\omega]_{\rm CE}.
\end{equation}

Recall now that the box constant is
\begin{equation}\label{30}
[W_\tau,\sigma_\omega]_{\rm Box}
=
\sup_{\gamma\in T_\omega^n}
\frac{
\displaystyle
\sum_{\alpha\leq\gamma}
W_\tau(\alpha)
(I^*\sigma_\omega)(\alpha)^2
}{
(I^*\sigma_\omega)(\gamma)
}
=
\sup_{\substack{\gamma\in T_\omega^n\\
\nu(S_\omega(\gamma))>0}}
\frac{
\displaystyle
\sum_{\alpha\leq\gamma}
W_\tau(\alpha)\nu(S_\omega(\alpha))^2
}{
\nu(S_\omega(\gamma))
}.
\end{equation}

We now verify the corner condition for $W_\tau$. Fix a split $S\subseteq\{1,\ldots,n\}$, and take vertices
$a\leq u$ in the $S$-coordinates and $b\leq v$ in the
complementary coordinates. Set
\[
s=
1+\sum_{j\in S}d_j(a)
+\sum_{j\notin S}d_j(b), 
\qquad
x=\sum_{j\in S}\bigl(d_j(u)-d_j(a)\bigr)
\quad
\text{and} \quad
y=\sum_{j\notin S}\bigl(d_j(v)-d_j(b)\bigr).
\]

Since $a\leq u$ and $b\leq v$, we have $x,y\geq0$. Thus,
\[
W_\tau(a,b)=s^{-\tau},\ 
W_\tau(u,b)=(s+x)^{-\tau},\ 
W_\tau(a,v)=(s+y)^{-\tau}
\text{ and } 
W_\tau(u,v)=(s+x+y)^{-\tau}.
\]

Therefore, since 
\(
(s+x)(s+y)-s(s+x+y)=xy\geq0,
\)
\[
\begin{aligned}
\frac{W_\tau(a,v)W_\tau(u,b)}
{W_\tau(a,b)W_\tau(u,v)}
&=
\left(
\frac{s(s+x+y)}{(s+x)(s+y)}
\right)^\tau \leq1,
\end{aligned}
\]

Thus $W_\tau$ satisfies the unit corner condition for every
split $S$. By Theorem \ref{thm:corner}, we obtain
\[
[W_\tau,\sigma_\omega]_{\rm Box}
\leq
[W_\tau,\sigma_\omega]_{\rm CE}
\leq
4^n[W_\tau,\sigma_\omega]_{\rm Box}.
\]

Taking the supremum over $\omega$, and using \eqref{27}, gives
\[
\|T_{K_\tau}\|
\approx_{n,\tau}
\sup_{\omega\in\mathbb T^n}
[W_\tau,\sigma_\omega]_{\rm Box}.
\]
which is what we wanted to prove.

Finally, by Theorem \ref{thm:capacity-characterization}, we also obtain
\[
\|T_{K_\tau}\|
\approx_{n,\tau}
\sup_{\omega\in\mathbb T^n}
C_{\rm cap}(W_\tau,\sigma_\omega).
\]

For a general finite measure $\nu$, set
$\nu_r=\mathbf 1_{r\mathbb D^n}\nu$.
Monotone convergence gives convergence of
the sharp embedding constants as
$r$ increases to $1$. Thus \eqref{27} extends to $\nu$, and
the corner and capacity theorems applied
to $\sigma_\omega$ complete the proof.

\end{proof}

\end{document}